\documentclass[a4paper,reqno,twoside,draft]{amsart}

\usepackage{latexsym}

\theoremstyle{plain}
\newtheorem{propn}{Proposition}[section]
\newtheorem{thm}[propn]{Theorem}
\newtheorem{lemma}[propn]{Lemma}

\newtheorem*{thm*}{Theorem}

\theoremstyle{definition}

\theoremstyle{remark}
\newtheorem*{rem}{Remark}
\newtheorem*{rems}{Remarks}

\newcommand{\Fock}{\mathcal{F}}
\newcommand{\Exp}{\mathcal{E}}
\newcommand{\e}[1]{\varepsilon (#1)}
\newcommand{\w}[1]{\varpi (#1)}
\newcommand{\step}{\mathbb{S}}

\newcommand{\init}{\mathfrak{h}}
\newcommand{\D}{\mathfrak{D}}

\newcommand{\noise}{\mathsf{k}}
\newcommand{\khat}{\wh{\noise}}

\newcommand{\iD}{\mathsf{D}}
\newcommand{\iDhat}{{\wh{\iD}}}
\newcommand{\Ihat}{\wh{I}}

\newcommand{\Q}{\mathcal{Q}}

\newcommand{\Qcd}{Q^{c,d}}
\newcommand{\Qab}{Q^{(\al,\be)}}

\newcommand{\Gab}{\Gh{\al}{\be}}

\newcommand{\Gh}[2]{G_{({#1},{#2})}}

\newcommand{\tot}{_{[0,t[}}

\newcommand{\fromt}{_{[t,\infty[}}

\newcommand{\ind}{\mathbf{1}}

\newcommand{\Op}{\mathcal{O}}
\newcommand{\Mat}{\mathrm{M}}
\newcommand{\bfc}{\mathbf{c}}

\newcommand{\Fab}{F^\al_\be}

\newcommand{\Lambdaab}{\Lambda^\al_\be}

\newcommand{\al}{\alpha}
\newcommand{\be}{\beta}

\newcommand{\si}{\sigma}
\newcommand{\lam}{\lambda}

\newcommand{\bvarpi}{\boldsymbol{\varpi}}

\newcommand{\Real}{\mathbb{R}}
\newcommand{\Rplus}{\Real_+}
\newcommand{\Comp}{\mathbb{C}}
\newcommand{\Nat}{\mathbb{N}}
\newcommand{\Int}{\mathbb{Z}}

\newcommand{\ip}[2]{\langle #1, #2 \rangle}

\newcommand{\norm}[1]{\lVert #1 \rVert}

\newenvironment{spmatrix}
{\bigl(\begin{smallmatrix}}{\end{smallmatrix}\bigr)}

\newcommand{\wh}{\widehat}
\newcommand{\wt}{\widetilde}
\newcommand{\ol}{\overline}

\newcommand{\schur}{\bullet}
\newcommand{\ot}{\otimes}

\newcommand{\op}{\oplus}
\newcommand{\uot}{\underline{\otimes}}

\newcommand{\To}{\rightarrow}
\newcommand{\Tends}{\rightarrow}
\newcommand{\Implies}{\Rightarrow}

\newcommand{\ti}{\textit}
\newcommand{\tu}{\textup}
\newcommand{\dfn}{\ti}

\DeclareMathOperator{\Dom}{Dom}
\DeclareMathOperator{\Ran}{Ran}
\DeclareMathOperator{\Lin}{Lin}

\DeclareMathOperator{\id}{id}
\DeclareMathOperator{\re}{Re}

\newenvironment{rlist}
{

\begin{enumerate}}
{\end{enumerate}}

\newcounter{step_count}
\newcommand{\Step}[1]{\bigskip
\noindent\stepcounter{step_count}
\arabic{step_count}. \ti{#1}: \ }

\numberwithin{equation}{section}

\begin{document}

\title[Quantum stochastic operator cocycles]
{Construction of some Quantum Stochastic Operator Cocycles
by the Semigroup Method}
\author{J.~Martin Lindsay}
\address{Department of Mathematics and Statistics \\ Lancaster
University \\ Lancaster LA1 4YF \\ UK}
\email{j.m.lindsay@lancaster.ac.uk}
\author{Stephen J.~Wills}
\address{School of Mathematical Sciences \\ University College Cork \\
Cork \\ Ireland}
\email{s.wills@ucc.ie}
\dedicatory{In celebration of Kalyan Sinha's sixtieth birthday}

\begin{abstract}
A new method for the construction of Fock-adapted operator Markovian
cocycles is outlined, and its use is illustrated by application to a
number of examples arising in physics and probability. The
construction uses the Trotter-Kato Theorem and a recent
characterisation of such cocycles in terms of an associated family of
contraction semigroups.
\end{abstract}

\maketitle

\section*{Introduction}

Let $\init$ and $\noise$ be fixed but arbitrary Hilbert spaces, and
let $\Fock$ denote the symmetric Fock space over $L^2 (\Rplus;
\noise)$. In this paper we are concerned with finding \dfn{left
contraction cocycles on $\init$ with noise dimension space $\noise$}.
That is families $V = (V_t)_{t \geq 0}$ of contractions on $\init \ot
\Fock$ that satisfy the functional equation
\[
V_0 = I, \qquad V_{r+t} = V_r \si_r (V_t) \quad \text{for all }
r,t \geq 0.
\]
Here $(\si_t)_{t \geq 0}$ is the endomorphism semigroup of right
shifts on the algebra $B(\init \ot \Fock)$ and the process $V$ is
required to be adapted to the natural filtration of subalgebras.

Subject to the further requirement that the cocycle be
\ti{Markov-regular}, there is a one-to-one correspondence between such
cocycles and solutions of a class of quantum stochastic differential
equations of Hudson-Parthasarathy type. The passage from cocycle to
solution of QSDE and back can often still be accomplished without the
restrictive Markov-regularity assumption, and so our results also
provide a new method for solving QSDEs with unbounded coefficients. In
common with earlier work of Fagnola and Mohari~(\cite{FFQP8, Moh}) the
resulting cocycle is realised as the limit (in the weak operator
topology) of a sequence of Markov-regular cocycles, obtained by
Yosida-type regularisation of the coefficients of the QSDE. These
authors approximate the cocycle directly, proving the existence of the
limit cocycle by appealing to the Arzel\`{a}-Ascoli Theorem. This
necessitates an assumption of separability on both of the Hilbert
spaces $\init$ and $\noise$. In contrast, our method goes via an
associated family of semigroups and employs the Trotter-Kato Theorem
for the approximations. No separability assumptions are needed here;
however there are some additional hypotheses. These amount to either
assuming relative boundedness of the annihilation coefficients with
respect to the time coefficient, or making operator core assumptions
on certain affine combinations of annihilation and time coefficients
of the QSDE. One goal of this paper is to show that such conditions
are natural and easy to check. Both our work and the cited papers of
Fagnola and Mohari circumvent limitations of earlier work on the
existence of solutions to QSDEs (\cite{FFPTRF,VS,App}). In order to
realise solutions by Picard iteration, as when the coefficients are
bounded, those papers imposed strong requirements such as the
existence of a dense subspace $\D$ on $\init$ that is left invariant
by all of the coefficients and on which these operators satisfy
suitable growth conditions.

\section{Cocycles and Semigroups}

In contrast with~\cite{theory}, which is coordinate-free, our account
in this paper will exploit a \ti{fixed orthonormal basis} $\eta =
\{d_i\}_{i \in I}$ of the noise dimension space $\noise$. We use an
index set $I$ not containing $0$ and use $\al, \be$ etc.\ to denote
elements of the augmented set $\Ihat := I \cup \{0\}$, reserving $i$
and $j$ for elements of the set $I$. Let
\begin{gather}
d_0 := 0, \quad \iD = \Lin \{d_i\}_{i \in I} \subset \noise \notag \\
\khat := \Comp \op \noise, \quad \iDhat = \Comp \op \iD, \quad \wh{c} =
\begin{spmatrix} 1 \\ c \end{spmatrix}, \label{hat defns} \\
e_0 = \wh{d_0}, \quad e_i = \begin{spmatrix} 0 \\ d \end{spmatrix}
\text{ where } d = d_i, \quad \Delta = I_\init \ot P_\noise, \notag
\end{gather}
for $c \in \noise$ and $i \in I$, and where $P_\noise \in B(\khat)$ is
the orthogonal projection with range $\noise$. In particular,
$\wh{\eta} := \{e_\al\}_{\al \in \Ihat}$ is thereby an orthonormal
basis of $\khat$ that contains the vector $\begin{spmatrix} 1 \\ 0
\end{spmatrix}$ and determines an isomorphism:
\[
\init \ot \khat \cong \bigoplus_{\al \in \Ihat} \init \quad \text{
via } \quad \xi \leftrightarrow (\xi^\al)_{\al \in \Ihat} \ \text{
where } \ \xi^\al = \big(I_\init \ot \langle e_\al |\big)\, \xi.
\]
Thus, for any subspace $\D$ of $\init$, the algebraic tensor product
$\D \uot \iDhat$ corresponds to those vectors with only finitely many
nonzero components all taken from $\D$.

Suppose that $[\Fab]_{\al,\be \in \Ihat}$ is a matrix of operators on
$\init$ and $\D$ is a subspace of $\init$ contained in the domain of
each $\Fab$. Then a sesquilinear form is defined on $\D \uot \iDhat$
through
\[
\bigl( (u^\al), (v^\be) \bigr) \mapsto \sum_{\al,\be} \ip{u^\al}{\Fab
v^\be}.
\]
Such a matrix can arise as the component matrix of an operator $F$ in
$\Op (\D \uot \iDhat)$, the linear space of operators on $\init \ot
\khat$ with domain $\D \uot \iDhat$, whose $(\al,\be)$-component is
determined by
\[
\ip{u}{\Fab v} = \ip{u \ot e_\al}{F v \ot e_\be}, \qquad u \in
\init, v \in \D, \al,\be \in \Ihat.
\]
If $\noise$ is finite dimensional then clearly all matrices of
operators with common dense domain $\D$ arise in this way. For
infinite dimensional $\noise$ the infinite matrices $[\Fab]$ that
correspond to an operator on $\D \uot \iDhat$ are precisely those that
satisfy
\[
\sum_\al \norm{\Fab u}^2 < \infty \quad \text{ for all } u \in \D, \be
\in \Ihat.
\]
We call such matrices \dfn{semiregular}, extending the terminology
introduced in~\cite{mother}. In the sequel we identify semiregular
matrices with their corresponding operators in $\Op (\D \uot \iDhat)$.

\begin{thm} \label{construct}
Let $\D$ be a dense subspace of $\init$. Let $F = [\Fab] \in \Op (\D
\uot \iDhat)$ satisfy
\begin{equation} \label{form ineq}
2 \re \ip{\xi}{F \xi} + \norm{\Delta F \xi}^2 \leq 0 \quad \text{for
all } \xi \in \D \uot \iDhat.
\end{equation}
Suppose that $F^0_0$ is the pregenerator of a $C_0$-semigroup and
either
\begin{rlist}
\item \label{rel bd}
$F^0_i$ is relatively bounded with respect to $F^0_0$ with relative
bound $0$ for each $i \in I$, or
\item \label{pregen}
$F^0_0 +F^0_i - \frac{1}{2}I_\init$ is the pregenerator of a
$C_0$-semigroup for each $i \in I$.
\end{rlist}
Then the operator QSDE
\begin{equation} \label{LHP}
dV_t = \sum_{\al,\be} V_t \Fab \, d\Lambdaab (t), \quad V_0 =
I_{\init \ot \Fock}
\end{equation}
has a strong solution $V$ on the domain $\D \uot \Exp_\eta$. Moreover
$V$ is a left contraction cocycle and is the unique contractive weak
solution on that domain.
\end{thm}

The following remarks explain the terminology used in the statement of
the theorem, and sketch details of the proof, a full account of which
can be found in~\cite{theory}.

\Step{Solutions of QSDEs}
Here $\Exp_\eta$ denotes the linear span of exponential vectors $\e{f}
= (1, f, (2!)^{-1/2} f \ot f, \ldots) \in \Fock$ where $f \in
\step_\eta$, the set of right continuous step functions on $\Rplus$,
with compact support, taking values in $\eta \cup \{0\}$. This is
dense in the Fock space $\Fock$ (\cite{L}, Proposition~2.1).

Three types of QSDE/solution are considered in~\cite{theory}, namely
the \dfn{form QSDE}, and \dfn{weak} or \dfn{strong solutions} of the
\dfn{operator QSDE}. The distinction between form QSDE and operator
QSDE rests on whether or not the coefficient matrix $[\Fab]$ is
semiregular, and it can be shown (\cite{theory},~Theorem~2.2) that
semiregularity of $[\Fab]$ is a necessary condition for the existence
of \ti{contractive} solutions of the form version of~\eqref{LHP}. For
contraction processes (and the \emph{left} Hudson-Parthasarathy QSDE
considered here), the distinction between weak and strong solutions
rests on whether or not the maps $t \mapsto V_t \xi$ are weakly or
strongly measurable. Our existence result yields a weakly continuous
contraction cocycle which is shown to be strongly continuous
(\cite{theory}, Lemma 1.2; \cite{jtl}, Proposition~1.1), and thus
strongly measurable.

Solutions of~\eqref{LHP} are characterised by their matrix elements
through
\begin{align*}
\ip{u \ot \e{f}}{(V_t-I) v \ot \e{g}} &= \int^t_0 \ip{u \ot \wh{f}(s)
\ot \e{f}}{\wh{V}_s (F \uot I_\Fock) v \ot \wh{g}(s) \ot \e{g}} \, ds
\\
&= \sum_{\al,\be} \int^t_0 \overline{f^\al (s)} g^\be (s) \ip{u \ot
\e{f}}{V_s \Fab v \ot \e{g}} \, ds
\end{align*}
for $u \in \init, v \in \D, f,g, \in \step_\eta$, where $\wh{f}(s) :=
\wh{f(s)} \in \khat$ is defined through~\eqref{hat defns}, $f^\al (s)
:= \ip{e_\al}{f(s)}$, a component of $\wh{f} (s)$ with respect to the
basis $\wh{\eta}$, and $\wh{V}_s$ is $V_s$ ampliated to $\init \ot
\khat \ot \Fock$.

\Step{Cocycles and semigroups}
There are two key results connecting cocycles to families of
semigroups. The first (from~\cite{father}) uses semigroups to
determine which contraction processes are cocycles; the second
(from~\cite{spawn}) characterises those families of semigroups that
can arise in this way, and so can be used to (re)construct a cocycle.

First, let $V = (V_t)_{t \geq 0}$ be an \dfn{adapted} process of
contraction operators, that is $V_t \in B(\init \ot \Fock\tot) \ot
I\fromt$ where $\Fock_J$ denotes the symmetric Fock space over $L^2
(J; \noise)$ and $I_J$ is the identity operator on $\Fock_J$. For any
$c,d \in \noise$ define contractions $\Qcd_t \in B(\init)$ through
\begin{equation} \label{semi defn}
\ip{u}{\Qcd_t v} = \ip{u \ot \w{c \ind\tot}}{V_t v \ot \w{d
\ind\tot}},
\end{equation}
where $\w{f} = \norm{\e{f}}^{-1} \e{f}$, the normalised exponential
vector associated to $f$. From~\cite{father} it follows that the
process $V$ is a left contraction cocycle if and only if $(\Qcd_t)_{t
\geq 0}$ is a semigroup on $\init$ for each pair $c,d \in \noise$, and
\begin{equation} \label{semi decomp}
\ip{u \ot \w{f \ind\tot}}{V_t v \ot \w{g \ind\tot}} =
\ip{u}{Q^{f(t_0), g(t_0)}_{t_1-t_0} \cdots Q^{f(t_n), g(t_n)}_{t-t_n}
v}
\end{equation}
for all $u,v \in \init$ and $f,g \in \step_\eta$, where $\{0 = t_0
\leq t_1 \leq \cdots \leq t_n \leq t\}$ contains the discontinuities
of $f \ind\tot$ and $g \ind\tot$.

This characterisation is independent of continuity assumptions in $t$.
However weak continuity of the map $t \mapsto V_t$ is equivalent to
strong continuity, which is equivalent to strong continuity of
\emph{all} (equivalently \emph{any}) of the semigroups $\Qcd$. In
particular if we assume this continuity and set $\Qab = \Qcd$ for $c =
d_\al, d = d_\be$, then each of these $C_0$-semigroups has a generator
$\Gab$. If $\D := \bigcap_{\al,\be} \Dom \Gab$ is a dense subspace of
$\init$ then it is not hard to show that $V$ satisfies the
QSDE~\eqref{LHP} on $\D \uot \Exp_\eta$ where the $\Fab$ are given by
\begin{align}
F^0_0 &= \Gh{0}{0}, \notag \\
F^i_0 &= \Gh{i}{0} -\Gh{0}{0} +\tfrac{1}{2}I_\init, \qquad F^0_j =
\Gh{0}{j}
-\Gh{0}{0} +\tfrac{1}{2}I_\init, \label{G defns} \\
F^i_j &= \Gh{i}{j} - \Gh{i}{0} -\Gh{0}{j} +\Gh{0}{0} -\delta^i_j
\notag
\end{align}
for $i,j \in I$, and where $\delta^i_j$ is the Kronecker delta. This 
is Theorem~4.1 of~\cite{theory}, proved there using a more direct 
method for extracting the $\Fab$. It extends Theorem~6.7 
of~\cite{father} when one, and hence all, of the semigroups is 
assumed to be norm continuous.

Thus associated to any cocycle is a family of semigroups. Conversely,
the following characterisation is obtained in~\cite{spawn}. A family
$\Q = \{\Qcd: c,d \in \iD\}$ of semigroups on $\init$ is amongst those
associated to a left contraction cocycle through~\eqref{semi defn} if
and only if for all $n \in \Nat$, $Y \in \Mat_n (|\init\rangle)$ and
positive invertible matrices $A,B \in \Mat_n (\Comp)$,
\begin{multline} \label{recon ineqs}
\norm{A^{-1/2} Y B^{-1/2}} \leq 1 \ \text{ implies } \\
\norm{(A \schur \bvarpi^\bfc_t)^{-1/2} (Q^\bfc_t \schur Y) (B
\schur \bvarpi^\bfc_t)^{-1/2}} \leq 1 \ \text{ for all } \bfc \in
\iD^n, t \geq 0.
\end{multline}
Here $|\init\rangle = B(\Comp; \init)$, the column operator space
associated to $\init$ (so that $\Mat_n (|\init\rangle)$ is identified
with $B(\Comp^n; \init^n)$), $\bvarpi^\bfc_t$ is the scalar matrix
$[\ip{\w{c_i \ind\tot}} {\w{c_j \ind\tot}}] \in \Mat_n (\Comp)$,
$Q^\bfc_t$ is the operator matrix $[Q^{c_i,c_j}_t] \in B(\init^n)$,
and $\schur$ denotes the Schur product.

The usefulness of this result may be explained as follows. Suppose
that $(V^{(n)})_{n \geq 1}$ is a sequence of cocycles whose associated
semigroups converge in the strong operator topology to a new family
$\Q = \{\Qcd\}$ of semigroups (${}^{(n)} \Qcd_t u \Tends \Qcd_t u$ for
all $c,d \in \iD$, $t \geq 0, u \in \init$). Then, since~\eqref{recon
ineqs} is manifestly stable under pointwise limits, there must be a
left contraction cocycle whose associated semigroups include $\Q$.

\Step{Trotter-Kato Theorem for cocycles/Proof of
\tu{Theorem~\ref{construct}}}
Requiring contractivity of solutions of~\eqref{LHP} imposes further
useful necessary conditions on the putative `stochastic generator' $F
= [\Fab] \in \Op (\D \uot \iD)$, namely:
\begin{itemize}
\item
$F$ satisfies the form inequality~\eqref{form ineq};
\item
each $F^i_0$ is relatively bounded with respect to $F^0_0$ with
relative bound $0$;
\item
each $\Gab$ defined through~\eqref{G defns} is dissipative;
\item
$[F^i_j +\delta^i_j I]_{i,j \in I}$ defines a contraction on $\init
\ot \noise$;
\item
\ti{If} $F^0_0$ is a pregenerator of a $C_0$-contraction semigroup
then there is a sequence $(F^{(n)})_{n \geq 1}$ in $B(\init \ot
\khat)$ with each $F^{(n)}$ satisfying~\eqref{form ineq}, and such
that $F^{(n)} \Tends F$ strongly on $\D \uot \iD$.
\end{itemize}
Most of these are noted in~\cite{FFQP8}; they are all contained in
Propositions~2.1,~3.1 and~3.2 of~\cite{theory}. For each $n$, the
QSDE~\eqref{LHP} with coefficient $F^{(n)}$ can be solved by Picard
iteration, and since $F^{(n)}$ satisfies~\eqref{form ineq} the
solution $V^{(n)}$ is a Markov-regular left contraction cocycle. Thus
each $V^{(n)}$ has an associated family of semigroups $\Q^{(n)} =
\{{}^{(n)} \! \Qab\}$, whose corresponding generators $\Gab^{(n)}$ are 
affine
combinations of the components of $F^{(n)}$ (as in~\eqref{G defns}),
and so converge strongly to $\Gab$ on $\D$. In turn these are
\ti{pregenerators} of $C_0$-contraction semigroups --- this is deduced
from a combination of the necessary conditions on $F$ listed above and
either condition~\eqref{rel bd} or~\eqref{pregen} of
Theorem~\ref{construct}. The resulting family $\Q$
satisfies~\eqref{recon ineqs}, and so we obtain the desired left
contraction cocycle $V$; moreover $V^{(n)}_t \Tends V_t$ for each $t
\geq 0$ in the weak operator topology.

\bigskip

A different type of semigroup appears in the characterisation of the
coefficients $F$ which yield \ti{isometric} solutions for the
QSDE~\eqref{LHP}. If $F \in B(\init \ot \khat)$ then $V$ is isometric
if and only if $F+F^* +F^* \Delta F = 0$, that is, equality holds
in~\eqref{form ineq} (\cite{FFQP8,mother}). For unbounded generators
this condition is still necessary, but it is no longer sufficient. If
$V = (V_t)_{t \geq 0}$ is a contractive solution of~\eqref{LHP}, then
defining
\begin{equation} \label{QDS}
\mathcal{T}_t (X) := \mathbb{E} [V^*_t (X \ot I_\Fock) V_t], \qquad X
\in B(\init),
\end{equation}
where $\mathbb{E}$ is the vacuum conditional expectation $B(\init \ot
\Fock) \To B(\init)$, produces a family of normal completely positive
contractions on $B(\init)$. Moreover $(\mathcal{T}_t)_{t \geq 0}$ is a
\dfn{quantum dynamical semigroup}, so that $\mathcal{T}_{0} = 
\id_{B(\init)}$, $\mathcal{T}_{r+t} = \mathcal{T}_r \circ 
\mathcal{T}_t$, and $t \mapsto \mathcal{T}_t (X)$ is ultraweakly 
continuous. It is the minimal quantum dynamical semigroup with 
form-generator $\mathcal{L}$, where
\[
\mathcal{L} (X): (u,v) \mapsto \ip{u}{X F^0_0 v} +\ip{F^0_0 u}{Xv}
+\sum_i \ip{F^i_0 u}{X F^i_0 v}, \quad u,v \in \D,
\]
meaning that
\begin{equation} \label{int id}
\ip{u}{\bigl( \mathcal{T}_t (X)-X \bigr)v} = \int^t_0 \mathcal{L}
\bigl(\mathcal{T}_s (X) \bigr) (u,v) \, ds.
\end{equation}
and
\[
\mathcal{T}_t (X) \leq \mathcal{T}'_t(X) \quad \text{for all } t
\geq 0 \text{ and } X \in B(\init)_+
\]
for any other quantum dynamical semigroup $\mathcal{T}'$ that
satisfies~\eqref{int id}. See~\cite{Proy} for more details.

\begin{thm}[\cite{FFQP8,Proy}] \label{iso coc}
Suppose that $F \in \Op (\D \uot \iD)$ satisfies
\begin{equation} \label{iso gen equal}
2 \re \ip{\xi}{F \xi} +\norm{\Delta F \xi}^2 = 0 \qquad \text{for all
} \xi \in \D \uot \iD
\end{equation}
and that $F^0_0$ is a pregenerator of a $C_0$-contraction semigroup.
Suppose also that there is a contractive solution $V$ to~\eqref{LHP}
for this $F$, and define $\mathcal{T}$ by~\eqref{QDS}. Then $V$ is
isometric if and only if $\mathcal{T}$ is conservative, that is,
$\mathcal{T}_t (I) = I$ for all $t \geq 0$.
\end{thm}

A number of necessary and sufficient conditions have been found for
conservativity of the quantum dynamical semigroup $\mathcal{T}$;
however in practice these can be difficult to verify, so a number of
sufficient conditions have been developed (see~\cite{suff2} and
references therein).

To determine if $V$ is coisometric involves analysing its
\dfn{Journ\'{e} dual}. Taking adjoints exchanges left cocycles and
right cocycles (satisfying $U_{r+t} = \sigma_r (U_t) U_r$); similarly
conjugating by time-reversal operators converts a left cocycle into a
right cocycle and vice versa. Thus conjugating the adjoint of a left
cocycle $V$ yields another left cocycle $\wt{V}$, and if $V$
satisfies~\eqref{LHP} for $F = [\Fab]$ then, at least formally,
$\wt{V}$ should satisfy~\eqref{LHP} for the adjoint matrix $F^\dagger
:= [{F^\beta_\alpha}^*]$, so that one can seek to apply
Theorem~\ref{iso coc} to this cocycle. This has been carried out
successfully in many instances (\cite{Proy}).

\section{Examples}

\subsection{A nonisometric cocycle as a limit of isometric cocycles}

Let $\init = l^2(\Int_+)$ with its standard orthonormal basis
$(e_n)_{n \geq 0}$, let $W \in B(\init)$ be the isometric right shift
$We_n = e_{n+1}$, and let $L$ be its Cayley transform:
$i(I+W)(I-W)^{-1}$. Thus $L$ is closed, densely defined and symmetric,
but not self-adjoint. In view of the identity $\Dom L^* = \Dom L
+\Comp e_0$, if we set $\D = \Dom L^*L$ and $\noise = \Comp$, then
\[
\begin{bmatrix} -\frac{1}{2} L^*L & -L^* \\ L & 0
\end{bmatrix}
\]
restricts to an operator $F\in\Op (\D\op \D)$, and moreover $F$
satisfies~\eqref{iso gen equal}. Since $L^*L$ is a positive
self-adjoint operator, $F^0_0$ is the generator of a $C_0$-contraction
semigroup. Also, since $-L^* \supset -L$, $F^0_1$ is relatively
bounded with respect to $F^0_0$ with relative bound $0$. Thus
condition~\eqref{rel bd} of Theorem~\ref{construct} applies. Let $V$
be the solution to~\eqref{LHP} for this $F$.

The related quantum dynamical semigroup $\mathcal{T}$ is \ti{not}
conservative, as noted in Example~3.5 of~\cite{noncon}, hence $V$
\ti{cannot} be isometric by Theorem~\ref{iso coc}. However, $V$ is
constructed by regularising $F$ to produce a sequence $(F^{(n)})_{n
\geq 1}$ with $F^{(n)} \Tends F$ strongly, so that for the semigroup
generators $\Gab^{(n)} \Tends \Gab$ strongly on $\D$, hence ${}^{(n)}
\!\Qab_t \Tends \Qab_t$ strongly as well, by the Trotter-Kato Theorem
(\cite{EngNag}, Theorem~III.4.8). Now for each $n \in \Nat$, $F^{(n)}$
lies in $B(\init \op \init)$ and satisfies~\eqref{iso gen equal}, hence
each cocycle $V^{(n)}$ is isometric. Thus, although~\eqref{semi
decomp} readily guarantees the weak operator convergence $V^{(n)}_t
\Tends V_t$ in our stochastic generalisation of the Trotter-Kato
Theorem, this example shows that convergence of the cocycles in the 
strong operator topology cannot be guaranteed.

\subsection*{Technical interlude}

Let $\init = l^2(\Int)$, and let $W$ be the right shift operator.
Given functions $\lam, \mu:\Int \To \Comp$ we will consider operators
such as $\lam(N) +W\mu(N)$ and $\lam(N) +\mu(N)W^*$, where $N$ denotes
the number operator. If $(e_n)_{n \in \Int}$ denotes the standard
orthonormal basis of $\init$ then these operators leave $\D := \Lin
\{e_n\}$ invariant and their action is given there by $\bigl( \lam(N)
+W\mu(N) \bigr) e_n = \lam(n) e_n +\mu(n) e_{n+1}$. However, despite
the temptation to restrict everything to $\D$ and thus simplify
various analytical domain considerations, we shall also have need to
think of these operators acting on their maximal domains. For example:
$\lam(N) +W\mu(N)$ has maximal domain
\[
\Bigl\{u = (u_n) \in \init: \sum_n |\lam(n) u_n +\mu(n-1)
u_{n-1}|^2 < \infty \Bigr\}.
\]
Let $\lam(N) \dot{+} W\mu(N)$ denote the resulting operator. An
advantage of passing to these domains is that we then have the true
operator equalities
\[
\bigl(\lam(N) +W\mu(N)|_\D \bigr)^* = \bigl( \lam(N) \dot{+} W\mu(N)
\bigr)^* = \ol{\lam}(N) \dot{+} \ol{\mu}(N) W^*,
\]
the subspace $\D$ being a core for such operators.

\begin{rem}
In general the maximal domain is strictly larger than that of the
operator sum, namely $\Dom \lam(N) \cap \Dom \mu(N)$. Indeed, if we
take $\lam(n) = n$, $\mu(n) = -n$ and set $u = (u_n)$ where $u_n =
(1+|n|)^{-1}$ then $u \in \init$, $u \notin \Dom \lam(N) = \Dom
\mu(N)$, but
\[
\sum_{n \in \Int} |\lam(n) u_n +\mu(n-1) u_{n-1}|^2 = \sum_{n<0}
\bigl[ (1-n)(2-n) \bigr]^{-2} + \frac{1}{4} +\sum_{n>0} \bigl[
n(n+1) \bigr]^{-2} < \infty,
\]
so that $u$ lies in the maximal domain of $\lambda(N)+W\mu(N)$.
\end{rem}

\subsection{The inverse harmonic oscillator}

Let $\init = l^2 (\Int_+)$ with its standard orthonormal basis
$(e_n)_{n \geq 0}$, let $\D = \Lin \{e_n\}$, and again let $W$ be the
isometric right shift on $\init$. Define $F \in \Op (\D \op \D)$ by
\begin{equation} \label{IHO}
F = \begin{bmatrix} -\frac{1}{2} |\lam|^2 (N+1) +i\mu(N) & W^*
\ol{\lam} (N) \\ -\lam(N) W & 0 \end{bmatrix},
\end{equation}
where $\lam:\Int_+ \To \Comp$ and $\mu:\Int_+ \To \Real$ are arbitrary
functions. The components of $F$, viewed as operators on $\D$, each
leave $\D$ invariant and it is not hard to verify that $F$
satisfies~\eqref{iso gen equal}.

Now $F^0_0 \in \Op (\D)$ is dissipative, and $(\nu -F^0_0) (\D) = \D$
for each $\nu > 0$, hence $\ol{F^0_0}$ is the generator of a
$C_0$-contraction semigroup. If there is $c > 0$ such that
\begin{equation} \label{growth}
c |\lam(n)| \leq |\lam(n+1)| \qquad \text{for all } n \geq 0
\end{equation}
then, for each $u \in \D$ and $a > 0$, we have
\[
\norm{F^0_1 u}^2 \leq \ip{u}{|\lambda|^2(N)u} \leq \frac{a^2}{2}
\norm{u}^2 +\frac{1}{2a^2} \norm{\lam(N) u}^2 \leq \frac{a^2}{2}
\norm{u}^2 +\frac{2}{a^2 c^4} \norm{F^0_0 u}^2,
\]
and so condition~\eqref{rel bd} of the theorem is satisfied. A
particular case where~\eqref{growth} holds is if we take $\lam(n) =
\sqrt{n}$; this was considered in~\cite{Wal}. However examples where
condition~\eqref{rel bd} of the theorem is violated are easily
generated. For example, setting
\[
\mu = 0 \ \text{ and } \ \lam(n) = \begin{cases} \sqrt{n} & \text{if
$n$ is odd}, \\ 0 & \text{if $n$ is even}, \end{cases}
\]
then $\norm{F^0_1 e_{2m+1}} = \sqrt{2m+1}$ but $\norm{F^0_0 e_{2m+1}}
= 0$ for each $m \geq 0$.

In fact every $F$ of the form~\eqref{IHO} stochastically generates a
cocycle since condition~\eqref{pregen} of Theorem~\ref{construct}
always holds. This can be seen using the following special case of a
well-known result (\cite{EngNag}, Corollary~II.3.17), whose short 
proof we include for the convenience of the reader.

\begin{lemma} \label{adj gen}
Let $A$ be a dissipative operator with dense domain $\D$. If $\D 
\subset \Dom A^*$ and $\D$ is  core for $A^*$ then $A$ is a 
pregenerator of a $C_0$-contraction semigroup.
\end{lemma}

\begin{proof}
By the Lumer-Phillips Theorem it suffices to show that $I-A$ has dense
range. Let $v \in \Ran (I-A)^\perp$, then
\[
\ip{(I-A)u}{v} = 0 \ \text{ for all } u \in \D \ \Implies \ v \in
\Dom A^* \ \text{ with } \ (I-A^*)v = 0.
\]
Dissipitavity of $A$ and the core assumption imply that $A^*$ is also
dissipative, so
\[
\norm{(I-A^*)v}\norm{v} \geq \re \ip{(I-A^*)v}{v} \geq \norm{v}^2,
\]
hence $v = 0$ as required.
\end{proof}

The putative stochastic generator $F$ above satisfies the
equality~\eqref{iso gen equal}. It follows that the operator
\[
A := F^0_0 +F^0_1 -\tfrac{1}{2}I|_\D
\]
is dissipative. Moreover $\D$ is a core for $A^* = \frac{1}{2}
|\lam|^2 (N+1) \dot{+} i\mu (N) \dot{-} W^* \ol{\lam} (N) \dot{-}
\frac{1}{2}I$, from the technical remarks above, and so Lemma~\ref{adj
gen} applies. Thus condition~\eqref{pregen} of Theorem~\ref{construct}
holds.

\subsection{Birth and death processes}

Now let $\init = l^2 (\Int)$, again write the standard orthonormal
basis as $(e_n)_{n \in \Int}$, set $\D = \Lin \{e_n\}$ and let $W$ be
the unitary right shift. This time let $\noise = \Comp^2$ and let $F
\in \Op (\D \uot \Comp^3)$ be given by
\[
F = \begin{bmatrix} -\frac{1}{2} |\lam|^2 (N) -\frac{1}{2} |\mu|^2 (N)
& \ol{\lam} (N) W^* & \ol{\mu} (N) W \\ -\lam (N) & W^*-I & 0 \\ -\mu
(N) & 0 & W-I \end{bmatrix}
\]
where $\lam, \mu:\Int \To \Comp$ are arbitrary functions. Again the
structure of $F$ and the fact that each component of $F$ leaves $\D$
invariant both facilitate the verification of~\eqref{iso gen equal}.
Careful choice of $\mu$ and $\lam$ again shows that
condition~\eqref{rel bd} of Theorem~\ref{construct} need not hold.
However Lemma~\ref{adj gen} and the technical remarks once more imply
that $F^0_0 +F^0_1 -\tfrac{1}{2}I|_\D$ and $F^0_0 +F^0_2
-\tfrac{1}{2}I|_\D$ are indeed pregenerators of $C_0$-semigroups, so
that condition~\eqref{pregen} of Theorem~\ref{construct} is satisfied.

\begin{rem}
In general it is not always possible to apply Theorem~\ref{iso coc}
(to $F$ and its adjoint) to conclude that $V$ is unitary. The
functions $\lam$ and $\mu$ that lead to unitary solutions have been
studied in detail in~\cite{FFBandD}. Necessary growth conditions on
$\lambda$ and $\mu$ for unitarity of $V$ are obtained from a
reformulation of the well-known resolvent conditions of Feller for the
Kolmogorov equations in the classical theory of Markov processes ---
see~\cite{MoFest} for further remarks in this vein.
\end{rem}

\subsection{A system for second harmonic generation}

In~\cite{state diff} an example is given involving two oscillators
coupled by a nonlinear interaction term. To obtain a quantum
stochastic dilation take $\init = l^2 (\Int_+^2) = l^2 (\Int_+) \ot
l^2 (\Int_+)$ with orthonormal basis 
$\{e_{m,n} = e_m\ot e_n\}_{m,n \in \Int_+}$,
and define creation and annihilation operators on each copy of $l^2
(\Int_+)$, so that, for example,
\[
a_1^* e_{m,n} = \sqrt{m+1} e_{m+1,n}, \qquad a_2 e_{m,n} =
\begin{cases} 0 & \text{if } n=0, \\ \sqrt{n} e_{m,n-1} & \text{if } n
\geq 1. \end{cases}
\]
Let $\D = \Lin \{e_{m,n}\}$, $\noise = \Comp^2$ and define $F \in \Op
(\D \uot \Comp^3)$ by
\[
F = \begin{bmatrix} K & -a_1^* & -a_2^* \\ a_1 & 0 & 0
\\ a_2 & 0 & 0 \end{bmatrix}
\]
where
\[
K = \tfrac{1}{2} (a_1^* a_1 + a_2^* a_2) + \omega (a_1^* -a_1) +\lam
({a_1^*}^2 a_2 -a_1^2 a_2^*)
\]
for some $\omega, \lam \in \Real$. Now clearly $F$ leaves $\D \uot
\Comp^3$ invariant, and it is easy to verify that $F$
satisfies~\eqref{iso gen equal}. In particular $K \in \Op (\D)$ is
dissipative, $\D$ is a core for $K^*$, and so $\ol{K}$ is the
generator of a contraction semigroup by Lemma~\ref{adj gen}. By the
same reasoning so are $\ol{K -a_1^* -\frac{1}{2}I}$ and $\ol{K -a_2^*
-\frac{1}{2}I}$. Hence condition~\eqref{pregen} of
Theorem~\ref{construct} holds, and so $F$ generates a contraction
cocycle $V$.

\begin{rems}
(i) A more detailed analysis of $K$, or rather of $K - \omega (a_1^*
-a_1)$, as has been done by Fagnola (\cite{2nd harm}) using the
results of~\cite{suff2}, shows that the associated quantum dynamical
semigroup is conservative, and hence that the cocycle $V$ is isometric
by Theorem~\ref{iso coc}. Using the results of~\cite{InvState},
Fagnola has also established the existence of an invariant state for
the semigroup.

(ii) The stochastic generator $F$ may be perturbed by introducing
nonzero number/exchange coefficients, without changing the quantum
dynamical semigroup, or complicating the application of
Theorem~\ref{construct}, provided that the annihilation coefficients
are appropriately modified. For example, the above argument works
equally well if we take
\[
F = \begin{bmatrix} K & -a_1^* V_1 & -a_2^* V_2 \\ a_1 & V_1-I & 0 \\
a_2 & 0 & V_2-I \end{bmatrix}
\]
for unitary operators $V_1$ and $V_2$ which leave $\iD$ invariant, for
example operators permuting the basis vectors $\{e_{m,n}\}$.
\end{rems}

\medskip
\noindent \ti{ACKNOWLEDGEMENTS}. We are grateful to Franco Fagnola
for drawing our attention to~\cite{state diff}, and sharing his
results in this area. SJW acknowledges support from EU TMR Networks
HPRN-CT-2002-00279 and HPRN-CT-2002-00280.

\end{document}